# FINITE SAMPLE PROPERTIES OF MULTIPLE IMPUTATION ESTIMATORS

By Jae Kwang Kim[1]

*Yonsei University*

Finite sample properties of multiple imputation estimators under the linear regression model are studied. The exact bias of the multiple imputation variance estimator is presented. A method of reducing the bias is presented and simulation is used to make comparisons. We also show that the suggested method can be used for a general class of linear estimators.

**1. Introduction.** Multiple imputation, proposed by Rubin (1978), is a procedure for handling missing data that allows the data analyst to use standard techniques of analysis designed for complete data, while providing a method to estimate the uncertainty due to the missing data. Repeated imputations are drawn from the posterior predictive distribution of the missing values under the specified model given a suitable prior distribution.

Schenker and Welsh [(1988), hereafter SW] studied the asymptotic properties of multiple imputation in the linear-model framework, where the scalar outcome variable $Y_i$ is assumed to follow the model

$$Y_i = \mathbf{x}'_i \boldsymbol{\beta} + e_i,$$
$$e_i \overset{\text{i.i.d.}}{\sim} N(0, \sigma^2). \quad (1)$$

The $p$-dimensional $\mathbf{x}_i$'s are observed on the complete sample and are assumed to be fixed.

To describe the imputation procedure, we adopt matrix notation. Without loss of generality, we assume that the first $r$ units are the respondents. Let $\mathbf{y}_r = (Y_1, Y_2, \ldots, Y_r)'$ and $X_r = (\mathbf{x}_1, \mathbf{x}_2, \ldots, \mathbf{x}_r)'$. Also, let $\mathbf{y}_{n-r} = (Y_{r+1}, Y_{r+2}, \ldots, Y_n)'$ and $X_{n-r} = (\mathbf{x}_{r+1}, \mathbf{x}_{r+2}, \ldots, \mathbf{x}_n)'$. The suggested method of multiple imputation for model (1) is as follows:

---

Received September 2001; revised January 2003.

[1]Supported by Korea Research Foundation Grant KRF-2003-003-C00032.

*AMS 2000 subject classifications.* Primary 62D05; secondary 62J99.

*Key words and phrases.* Missing data, nonresponse, variance estimation.

[M1] For each repetition of the imputation, $k = 1, \ldots, M$, draw

$$\sigma_{(k)}^{*2} \mid \mathbf{y}_r \overset{\text{i.i.d.}}{\sim} (r-p)\hat{\sigma}_r^2 / \chi_{r-p}^2, \tag{2}$$

where $\hat{\sigma}_r^2 = (r-p)^{-1} \mathbf{y}_r'[I - X_r(X_r'X_r)^{-1}X_r']\mathbf{y}_r$.

[M2] Draw

$$\boldsymbol{\beta}_{(k)}^* \mid (\mathbf{y}_r, \sigma_{(k)}^*) \overset{\text{i.i.d.}}{\sim} N(\hat{\boldsymbol{\beta}}_r, (X_r'X_r)^{-1}\sigma_{(k)}^{*2}), \tag{3}$$

where $\hat{\boldsymbol{\beta}}_r = (X_r'X_r)^{-1}X_r'\mathbf{y}_r$.

[M3] For each missing unit $j = r+1, \ldots, n$ draw

$$e_{j(k)}^{**} \mid (\boldsymbol{\beta}_{(k)}^*, \sigma_{(k)}^*) \overset{\text{i.i.d.}}{\sim} N(0, \sigma_{(k)}^{*2}). \tag{4}$$

Then $Y_{j(k)}^{**} = \mathbf{x}_j \boldsymbol{\beta}_{(k)}^* + e_{j(k)}^{**}$ is the imputed value associated with unit $j$ for the $k$th imputation.

[M4] Repeat [M1]–[M3] independently $M$ times.

The above procedure assumes a constant prior for $(\boldsymbol{\beta}, \log \sigma)$ and an ignorable response mechanism in the sense of Rubin (1976).

At each repetition of the imputation, $k = 1, \ldots, M$, we can calculate the imputed version of the full sample estimators

$$\hat{\boldsymbol{\beta}}_{I(k),n} = \left( \sum_{i=1}^n \mathbf{x}_i \mathbf{x}_i' \right)^{-1} \left[ \sum_{i=1}^r \mathbf{x}_i y_i + \sum_{i=r+1}^n \mathbf{x}_i Y_{i(k)}^{**} \right]$$

and

$$\hat{V}_{I(k),n} = \left( \sum_{i=1}^n \mathbf{x}_i \mathbf{x}_i' \right)^{-1} \hat{\sigma}_{I(k),n}^2,$$

where

$$\hat{\sigma}_{I(k),n}^2 = (n-p)^{-1} \left[ \sum_{i=1}^r (Y_i - \mathbf{x}_i'\hat{\boldsymbol{\beta}}_{I(k),n})^2 + \sum_{i=r+1}^n (Y_{i(k)}^{**} - \mathbf{x}_i'\hat{\boldsymbol{\beta}}_{I(k),n})^2 \right].$$

The proposed point estimator for the regression coefficient based on $M$ repeated imputations is

$$\hat{\boldsymbol{\beta}}_{M,n} = M^{-1} \sum_{k=1}^M \hat{\boldsymbol{\beta}}_{I(k),n}. \tag{5}$$

The proposed estimator for the variance of the point estimator (5) is

$$\hat{V}_{M,n} = W_{M,n} + (1 + M^{-1})B_{M,n}, \tag{6}$$

where

$$W_{M,n} = M^{-1} \sum_{k=1}^M \hat{V}_{I(k),n} \tag{7}$$

and

(8) $$B_{M,n} = (M-1)^{-1} \sum_{k=1}^{M} (\hat{\boldsymbol{\beta}}_{I(k),n} - \hat{\boldsymbol{\beta}}_{M,n})(\hat{\boldsymbol{\beta}}_{I(k),n} - \hat{\boldsymbol{\beta}}_{M,n})'.$$

Rubin (1987) called $W_{M,n}$ the *within-imputation variance* and called $B_{M,n}$ the *between-imputation variance*. We call $\hat{V}_{M,n}$ Rubin's variance estimator.

SW studied the asymptotic properties of the point estimator (5) and its variance estimator (6). Under regularity conditions they showed that

(9) $$\lim_{n\to\infty} E(\hat{\boldsymbol{\beta}}_{M,n} - \boldsymbol{\beta}) = 0$$

and

(10) $$\lim_{n\to\infty} n\{E(\hat{V}_{M,n}) - \operatorname{Var}(\hat{\boldsymbol{\beta}}_{M,n})\} = 0,$$

where the reference distribution in (9) and (10) is the regression model (1) with an ignorable response mechanism.

Note that (9) and (10) require that the sample size $n$ go to infinity, for fixed $M$, $M > 1$. Finite sample properties are not discussed by SW. The next section gives finite sample properties of the multiple imputation estimators. In Section 3 a simple modified version of the SW method is proposed to minimize the finite sample bias of the multiple imputation variance estimator. In Section 4 extensions are made to a more general class of estimators. In Section 5 results of a simulation study are reported. In Section 6 concluding remarks are made.

**2. Finite sample properties.** The following lemma provides the covariance structure of the multiply-imputed data set generated by [M1]–[M4].

LEMMA 2.1. *Let $Y_i$ be the observed value of the $i$th unit, $i = 1, 2, \ldots, r$ ($r > p+2$), and let $Y_{j(k)}^{**}$ be the imputed value associated with the $j$th unit for the $k$th repetition of the multiple imputation generated by the steps* [M1]–[M4]. *Then, under model* (1) *with an ignorable response mechanism,*

(11) $$\operatorname{Cov}(Y_i, Y_{j(k)}^{**}) = \mathbf{x}_i'(X_r'X_r)^{-1}\mathbf{x}_j\sigma^2$$

*and*

(12) $$\operatorname{Cov}(Y_{i(k)}^{**}, Y_{j(s)}^{**}) = \begin{cases} (1+\lambda)\mathbf{x}_i'(X_r'X_r)^{-1}\mathbf{x}_i\sigma^2 + \lambda\sigma^2, & \text{if } i = j \text{ and } k = s, \\ (1+\lambda)\mathbf{x}_i'(X_r'X_r)^{-1}\mathbf{x}_j\sigma^2, & \text{if } i \neq j \text{ and } k = s, \\ \mathbf{x}_i'(X_r'X_r)^{-1}\mathbf{x}_j\sigma^2, & \text{if } k \neq s, \end{cases}$$

*where $\lambda = (r-p-2)^{-1}(r-p)$ and the expectations in* (11) *and* (12) *are taken over the joint distribution of model* (1) *and the imputation mechanism with the indices of respondents fixed.*

For the proof, see Appendix A.

Note that the imputed value $Y_{i(k)}^{**}$ can be decomposed into three independent components as

$$Y_{i(k)}^{**} = \mathbf{x}_i\hat{\boldsymbol{\beta}}_r + \mathbf{x}_i(\boldsymbol{\beta}_{(k)}^* - \hat{\boldsymbol{\beta}}_r) + e_{i(k)}^{**}. \tag{13}$$

The first component $\mathbf{x}_i\hat{\boldsymbol{\beta}}_r$ has mean $\mathbf{x}_i\boldsymbol{\beta}$ and variance $\mathbf{x}_i'(X_r'X_r)^{-1}\mathbf{x}_i\sigma^2$, the second component has mean zero and variance $\lambda\mathbf{x}_i'(X_r'X_r)^{-1}\mathbf{x}_i\sigma^2$ and the third component has mean zero and variance $\lambda\sigma^2$.

The following theorem gives the mean and variance of the point estimator $\hat{\boldsymbol{\beta}}_{M,n}$ of the regression coefficient and the mean of the multiple imputation variance estimator $\hat{V}_{M,n}$. Again, the expectations in the following theorem are taken over the joint distribution of model (1) and the imputation mechanism with the indices of respondents fixed.

THEOREM 2.1. *Under the assumptions of Lemma* 2.1,

$$E(\hat{\boldsymbol{\beta}}_{M,n}) = \boldsymbol{\beta}, \tag{14}$$

$$\text{Var}(\hat{\boldsymbol{\beta}}_{M,n}) = (X_r'X_r)^{-1}\sigma^2 + M^{-1}\lambda[(X_r'X_r)^{-1} - (X_n'X_n)^{-1}]\sigma^2, \tag{15}$$

$$E(W_{M,n}) = (X_n'X_n)^{-1}\{1 + (n-p)^{-1}(\lambda-1)(n-r)\}\sigma^2 \tag{16}$$

*and*

$$E(B_{M,n}) = \lambda[(X_r'X_r)^{-1} - (X_n'X_n)^{-1}]\sigma^2, \tag{17}$$

*where $\lambda = (r-p-2)^{-1}(r-p)$, $W_{M,n}$ is the within-imputation variability defined in (7) and $B_{M,n}$ is the between-imputation variability defined in (8). The bias of the multiple imputation variance estimator is*

$$\begin{aligned}E(\hat{V}_{M,n}) - \text{Var}(\hat{\boldsymbol{\beta}}_{M,n}) &= (X_n'X_n)^{-1}(n-p)^{-1}(\lambda-1)(n-r)\sigma^2 \\ &\quad + (\lambda-1)[(X_r'X_r)^{-1} - (X_n'X_n)^{-1}]\sigma^2.\end{aligned} \tag{18}$$

For the proof, see Appendix B.

As is observed from (15), the point estimator $\hat{\boldsymbol{\beta}}_{M,n}$ for infinite $M$ achieves the same efficiency of $\hat{\boldsymbol{\beta}}_r$, the estimator based on the respondents. In fact,

$$\begin{aligned}\lim_{M\to\infty}\hat{\boldsymbol{\beta}}_{M,n} &= \left(\sum_{i=1}^n \mathbf{x}_i\mathbf{x}'\right)^{-1}\left\{\sum_{i=1}^r \mathbf{x}_i[\mathbf{x}_i'\hat{\boldsymbol{\beta}}_r + (y_i - \mathbf{x}_i'\hat{\boldsymbol{\beta}}_r)]\right\} \\ &\quad + \left(\sum_{i=1}^n \mathbf{x}_i\mathbf{x}'\right)^{-1}\left\{\sum_{i=r+1}^n \mathbf{x}_i\left[\mathbf{x}_i'\hat{\boldsymbol{\beta}}_r + M^{-1}\sum_{k=1}^M (\mathbf{x}_i'(\boldsymbol{\beta}_{(k)}^* - \hat{\boldsymbol{\beta}}_r) + e_{i(k)}^{**})\right]\right\} \\ &= \hat{\boldsymbol{\beta}}_r,\end{aligned}$$

because $\sum_{i=1}^{r} \mathbf{x}_i(y_i - \mathbf{x}_i'\hat{\boldsymbol{\beta}}_r) = 0$ by standard regression theory, $\lim_{M \to \infty} M^{-1} \times \sum_{k=1}^{M}(\boldsymbol{\beta}_{(k)}^* - \hat{\boldsymbol{\beta}}_r) = 0$ by the law of large numbers and $\lim_{M \to \infty} M^{-1} \times \sum_{k=1}^{M} e_{i(k)}^{**} = 0$ by the law of large numbers. By (15) the variance of $\hat{\boldsymbol{\beta}}_{M,n}$ can be written as

$$(19) \quad \operatorname{Var}(\hat{\boldsymbol{\beta}}_{M,n}) = \operatorname{Var}(\hat{\boldsymbol{\beta}}_r) + M^{-1}\lambda[\operatorname{Var}(\hat{\boldsymbol{\beta}}_r) - \operatorname{Var}(\hat{\boldsymbol{\beta}}_n)].$$

The second part in the right-hand side of (19) is the increase in variance due to using $\hat{\boldsymbol{\beta}}_{M,n}$ instead of $\hat{\boldsymbol{\beta}}_r$. By (17) that increase can be unbiasedly estimated by $M^{-1}B_{M,n}$. Thus, an alternative estimator for the total variance of $\hat{\boldsymbol{\beta}}_{M,n}$ that is unbiased for (15) is

$$(20) \quad \hat{\operatorname{Var}}(\hat{\boldsymbol{\beta}}_r) + M^{-1}B_{M,n},$$

where $\hat{\operatorname{Var}}(\hat{\boldsymbol{\beta}}_r)$ is the standard variance estimator that treats the respondents as if they are the original sample.

In large samples, $\lambda \doteq 1$ and the bias of the multiple imputation variance estimator for the imputed regression coefficient is negligible. The bias term (18) is an exact bias for $r > p + 2$.

The total variance of $\hat{\boldsymbol{\beta}}_{M,n}$ can be decomposed into three parts:

$$\operatorname{Var}(\hat{\boldsymbol{\beta}}_{M,n}) = \operatorname{Var}(\hat{\boldsymbol{\beta}}_n) + \operatorname{Var}(\hat{\boldsymbol{\beta}}_r - \hat{\boldsymbol{\beta}}_n) + \operatorname{Var}(\hat{\boldsymbol{\beta}}_{M,n} - \hat{\boldsymbol{\beta}}_r).$$

The first part, the second part and the third part can be called the sampling variance, the variance due to missingness and the variance due to imputation, respectively. Rubin's multiple imputation uses $W_{M,n}$ to estimate the sampling variance, $B_{M,n}$ to estimate the variance due to missingness and $M^{-1}B_{M,n}$ to estimate the variance due to $M$ repeated imputations. The first term on the right-hand side of (18) is the bias of $W_{M,n}$ as an estimator of the sampling variance and the second term on the right-hand side of (18) is the bias of $B_{M,n}$ as an estimator of the variance due to missingness. The imputation variance is unbiasedly estimated by $M^{-1}B_{M,n}$.

**3. A modification.** We consider ways of reducing the bias of the multiple imputation variance estimator. Recall that multiple imputation is characterized by the method of generating the imputed values and by the variance formula. The variance formula directly uses the complete sample variance estimator so that it can be implemented easily using existing software.

One estimator of variance is the alternative variance estimator in (20), which involves direct estimation of parameters from the respondents. A similar idea was used by Wang and Robins (1998). Then Rubin's variance estimator $\hat{V}_{M,n}$ in (6) is no longer required.

If we want to use Rubin's variance formula, one approach is to modify the imputation method to minimize the bias term in (18). To find a best imputation procedure, instead of fixing a constant prior for $\log \sigma$, we use a class of prior distributions indexed by hyperparameters to express a class of imputation methods. As a conjugate prior distribution for $\sigma^2$, we choose a scaled inverse chi-square with degrees of freedom $\nu_0$ and scale parameter $\sigma_0^2$ as the hyperparameters; that is, the prior distribution of $\sigma^2$ is the distribution of $\nu_0 \sigma_0^2 / \chi_{\nu_0}^2$. In the modified imputation method we determine the values of the hyperparameters, $\nu_0$ and $\sigma_0^2$, to remove the bias of the multiple imputation variance estimator.

Using the hyperparameters, the posterior distribution for $\sigma^2$ is written

$$\sigma_{(k)}^{*2} \mid \mathbf{y}_r \overset{\text{i.i.d.}}{\sim} [\nu_0 \sigma_0^2 + (r-p)\hat{\sigma}_r^2]/\chi_{(\nu_0+r-p)}^2, \tag{21}$$

so that

$$E(\sigma_{(k)}^{*2}) = (\nu_0 + r - p - 2)^{-1}(\nu_0 \sigma_0^2 + (r-p)\sigma^2). \tag{22}$$

Note that SW used $\nu_0 = 0$. Using the arguments of the proofs of Lemma 2.1 and Theorem 2.1, the bias of the multiple imputation variance estimator based on the posterior distribution in (21) is

$$\begin{aligned}
\text{Bias}(\hat{V}_{M,n}) = {} & (X_n' X_n)^{-1}(n-p)^{-1}\{(\lambda_0 - 1)(n-r)\sigma^2 + \lambda_1(n-r)\sigma_0^2\} \\
& + \{(X_r' X_r)^{-1} - (X_n' X_n)^{-1}\}\{(\lambda_0 - 1)\sigma^2 + \lambda_1 \sigma_0^2\},
\end{aligned} \tag{23}$$

where $\lambda_0 = (\nu_0 + r - p - 2)^{-1}(r-p)$ and $\lambda_1 = (\nu_0 + r - p - 2)^{-1}\nu_0$. The bias is zero when $\nu_0 = 2$ and $\sigma_0^2 = 0$, which is equivalent to generating the posterior values of $\sigma^2$ from the inverse chi-square distribution with degrees of freedom $\nu = r - p + 2$ instead of $\nu = r - p$ in (2). Thus, the choice of $\nu = r - p + 2$ in (2) makes Rubin's variance estimator unbiased for a finite sample.

We can also derive the optimal prior using the recent work of Meng and Zaslavsky (2002) on single observation unbiased priors (SOUP). Meng and Zaslavsky [2002, Section 6] showed that

$$\pi(\sigma^2) \propto (\sigma^2)^{-2} \tag{24}$$

is the unique SOUP among all continuously relatively scale invariant priors for the scale family distribution. Note that the prior density of the scaled inverse chi-square distribution can be written as

$$\pi(\sigma^2) \propto (\sigma^2)^{-(\nu_o/2+1)} \exp[-\nu_0 \sigma_0^2/(2\sigma^2)]. \tag{25}$$

Thus, the unique SOUP for $\sigma^2$ in (24) corresponds to the scaled inverse chi-square distribution with $\nu_0 = 2$ and $\sigma_0^2 = 0$, which makes the posterior mean in (22) unbiased for $\sigma^2$.

**4. Extensions.** In this section we investigate the properties of the modified method when it is applied to estimators other than the regression coefficients. Let the complete sample point estimator be a linear estimator of the form

$$\hat{\theta}_n = \sum_{i=1}^{n} \alpha_i Y_i \tag{26}$$

for some known coefficients $\alpha_i$. Also, let the complete sample estimator for the variance of $\hat{\theta}_n$ be a quadratic function of the sample values of the form

$$\hat{V}_n = \sum_{i=1}^{n} \sum_{j=1}^{n} \Omega_{ij} Y_i Y_j \tag{27}$$

for some known coefficients $\Omega_{ij}$.

For the $\hat{V}_{M,n}$ to be asymptotically unbiased for the variance of $\hat{\theta}_{M,n}$, we need the "congeniality" assumption as defined in Meng (1994). The congeniality assumption in our context implies

$$\operatorname{Var}(\hat{\theta}_{\infty,n}) = \operatorname{Var}(\hat{\theta}_n) + \operatorname{Var}(\hat{\theta}_{\infty,n} - \hat{\theta}_n). \tag{28}$$

For the case of $\hat{\theta}_n = \hat{\boldsymbol{\beta}}_n$ discussed in Section 2, congeniality holds because $\hat{\theta}_{\infty,n} = \hat{\boldsymbol{\beta}}_r$ and

$$\operatorname{Var}(\hat{\boldsymbol{\beta}}_r - \hat{\boldsymbol{\beta}}_n) = (X_r' X_r)^{-1} \sigma^2 - (X_n' X_n)^{-1} \sigma^2$$
$$= \operatorname{Var}(\hat{\boldsymbol{\beta}}_r) - \operatorname{Var}(\hat{\boldsymbol{\beta}}_n).$$

We restrict our attention to the case of congenial multiple imputation estimation because otherwise the variance estimator will be biased even asymptotically.

The following theorem expresses the bias of Rubin's variance estimator applied to the general class of estimators $\hat{\theta}_n$ in (26) and $\hat{V}_n$ in (27) under the SW method.

THEOREM 4.1. *Let the assumptions of Lemma* 2.1 *hold. Assume that the complete sample variance estimator $\hat{V}_n$ in* (27) *is unbiased for the variance of the complete sample point estimator $\hat{\theta}_n$ in* (26) *under model* (1). *Let the congeniality assumption* (28) *hold. Then the multiple imputation point estimator $\hat{\theta}_{M,n}$ is unbiased with variance*

$$\operatorname{Var}(\hat{\theta}_{M,n}) = \left\{ \sum_{i=1}^{r} \alpha_i^2 + 2 \sum_{i=1}^{r} \sum_{j=r+1}^{n} \alpha_i \alpha_j h_{ij} + \sum_{i=r+1}^{n} \sum_{j=r+1}^{n} \alpha_i \alpha_j h_{ij} \right\} \sigma^2$$
$$+ M^{-1} \lambda \left\{ \sum_{i=r+1}^{n} \sum_{j=r+1}^{n} \alpha_i \alpha_j h_{ij} + \sum_{i=r+1}^{n} \alpha_i^2 \right\} \sigma^2. \tag{29}$$

*The bias of $\hat{V}_{M,n}$ as an estimator of $\mathrm{Var}(\hat{\theta}_{M,n})$ is*

$$
\begin{aligned}
(30)\quad E(\hat{V}_{M,n}) - \mathrm{Var}(\hat{\theta}_{M,n}) &= \bigg\{ 2\sum_{i=1}^{r}\sum_{j=r+1}^{n}\Omega_{ij}h_{ij} + (1+\lambda)\sum_{i=r+1}^{n}\sum_{j=r+1}^{n}\Omega_{ij}h_{ij} \\
&\quad + (\lambda-1)\sum_{j=r+1}^{n}\Omega_{jj}\bigg\}\sigma^2 \\
&\quad + (\lambda-1)\bigg\{\sum_{j=r+1}^{n}\alpha_j^2 + \sum_{j=r+1}^{n}\sum_{j=r+1}^{n}\alpha_i\alpha_j h_{ij}\bigg\}\sigma^2,
\end{aligned}
$$

*where $\lambda = (r-p-2)^{-1}(r-p)$ and $h_{ij} = \mathbf{x}_i'(X_r'X_r)^{-1}\mathbf{x}_j$.*

For the proof, see Appendix C.

The first term on the right-hand side of (30) is the bias of $W_{M,n}$ as an estimator of $\mathrm{Var}(\hat{\theta}_n)$ and the second term is the bias of $(1+M^{-1})B_{M,n}$ as an estimator of $\mathrm{Var}(\hat{\theta}_{M,n} - \hat{\theta}_n)$.

Note that the bias term in (30) can be written

$$
(31)\quad \mathrm{Bias}(\hat{V}_{M,n}) = 2\sum_{i=1}^{n}\sum_{j=r+1}^{n}\Omega_{ij}h_{ij}\sigma^2 + (\lambda-1)U\sigma^2,
$$

where

$$
\begin{aligned}
U &= \bigg\{\sum_{i=r+1}^{n}\sum_{j=r+1}^{n}(\Omega_{ij}+\alpha_i\alpha_j)h_{ij} + \sum_{j=r+1}^{n}(\Omega_{jj}+\alpha_j^2)\bigg\} \\
&= \mathrm{trace}\{(\Omega_{n-r}+\boldsymbol{\alpha}_{n-r}\boldsymbol{\alpha}_{n-r}')[X_{n-r}(X_r'X_r)^{-1}X_{n-r}' + I_{n-r}]\}
\end{aligned}
$$

with $\Omega_{n-r}$ the lower-right $(n-r)\times(n-r)$ partition of $\Omega = [\Omega_{ij}]$ and $\boldsymbol{\alpha}_{n-r} = (\alpha_{r+1},\ldots,\alpha_n)'$. By the nonnegative definiteness of $\Omega_{n-r}$, $\boldsymbol{\alpha}_{n-r}\boldsymbol{\alpha}_{n-r}'$, $X_{n-r}(X_r'X_r)^{-1}X_{n-r}'$ and $I_{n-r}$, the $U$ term is nonnegative. Thus, if we use the modified method suggested in Section 3 so that we have $\lambda = 1$, the variance term (29) decreases and the bias will be reduced to

$$
(32)\quad \mathrm{Bias}(\hat{V}_{M,n}) = 2\sum_{i=1}^{n}\sum_{j=r+1}^{n}\Omega_{ij}h_{ij}\sigma^2,
$$

which is always smaller than the original bias in (31) because $\lambda = (r-p-2)^{-1}(r-p) > 1$.

A sufficient condition for the bias term in (32) to be zero is

$$
(33)\quad \Omega = c[I_n - X_n(X_n'X_n)^{-1}X_n']
$$

for some constant $c > 0$. To show this, let $\delta_{ij}$ be the $(i,j)$th element of $I_n$. Then

$$\sum_{i=1}^n \Omega_{ij} h_{ij} = \sum_{i=1}^n c[\delta_{ij} - \mathbf{x}_i'(X_n'X_n)^{-1}\mathbf{x}_j][\mathbf{x}_i'(X_r'X_r)^{-1}\mathbf{x}_j]$$

$$= c\mathbf{x}_j'(X_r'X_r)^{-1}\mathbf{x}_j - c\mathbf{x}_j'(X_n'X_n)^{-1}\left[\sum_{i=1}^n \mathbf{x}_i\mathbf{x}_i'\right](X_r'X_r)^{-1}\mathbf{x}_j = 0$$

and the bias term in (32) equals zero, which alternatively justifies our assertion of zero bias in Section 3.

**5. Simulation study.** To see the effect of changing $\nu = r - p$ into $\nu = r - p + 2$, we performed a limited simulation. The simulation study can be described as a $2 \times 3 \times 2$ factorial design with $L = 50{,}000$ samples within each cell, where each sample is generated from

(34) $$Y_i = 2 + 4x_i + e_i,$$

where $x_i = 5 + 10(n+1)^{-1}i$ and $e_i$ are independently and identically distributed from the standard normal distribution. Thus, the population mean of $(X, Y)$ is $(10, 42)$. The factors are as follows:

1. factor A, method of multiple imputation—SW method ($\nu = r - 2$), new method ($\nu = r$);
2. factor B, response rate $(r/n)$—0.8, 0.6, 0.4;
3. factor C, sample size $(n)$—20, 200.

We used a uniform response mechanism and $M = 5$ repeated imputations.

Table 1 presents the mean, the variance and the percentage relative efficiency of the point estimators under the two imputation schemes. The percentage relative efficiency (PRE) is

$$\text{PRE} = [\text{Var}_L(\hat{\theta}_{\text{SW}})]^{-1} \text{Var}_L(\hat{\theta}_{\text{new}}) \times 100,$$

where the subscript $L$ denotes the distribution generated by the Monte Carlo simulation. Both imputation methods are unbiased for the two parameters and the Monte Carlo results are consistent with that property. The new procedure is slightly more efficient than the SW procedure and the efficiency is greater for lower response rates.

Table 2 presents the relative bias and $z$-statistics for the variance estimators. The relative bias of $\hat{V}$ as an estimator of the variance of $\hat{\theta}$ is calculated as $[\text{Var}_L(\hat{\theta})]^{-1}[E_L(\hat{V}) - \text{Var}_L(\hat{\theta})]$, and the $z$-statistic for testing $H_0 : E(\hat{V}) = \text{Var}(\hat{\theta})$ is

(35) $$z\text{-}statistic = \frac{L^{1/2}[E_L\hat{V} - \text{Var}_L(\hat{\theta})]}{\{E_L[\hat{V} - E_L(\hat{V}) + \text{Var}_L(\hat{\theta}) - (\hat{\theta} - E_L(\hat{\theta}))^2]^2\}^{1/2}}.$$

TABLE 1
*Mean, variance and the percentage relative efficiency (PRE) of the multiple imputation point estimators under the two different imputation schemes (50,000 samples)*

| Parameter | $n$ | $r/n$ | Mean | | Variance | | PRE |
|---|---|---|---|---|---|---|---|
| | | | SW | New | SW | New | (%) |
| Mean | 20 | 0.8 | 42.0 | 42.0 | 0.066533 | 0.066056 | 99.2 |
| | | 0.6 | 42.0 | 42.0 | 0.095951 | 0.094192 | 98.2 |
| | | 0.4 | 42.0 | 42.0 | 0.150521 | 0.142537 | 94.7 |
| | 200 | 0.8 | 42.0 | 42.0 | 0.006594 | 0.006583 | 99.8 |
| | | 0.6 | 42.0 | 42.0 | 0.009069 | 0.009057 | 99.9 |
| | | 0.4 | 42.0 | 42.0 | 0.014143 | 0.014090 | 99.6 |
| Slope | 20 | 0.8 | 4.0 | 4.0 | 0.008873 | 0.008785 | 99.0 |
| | | 0.6 | 4.0 | 4.0 | 0.013148 | 0.012956 | 98.5 |
| | | 0.4 | 4.0 | 4.0 | 0.018190 | 0.017443 | 95.9 |
| | 200 | 0.8 | 4.0 | 4.0 | 0.000780 | 0.000779 | 99.9 |
| | | 0.6 | 4.0 | 4.0 | 0.001084 | 0.001084 | 100.0 |
| | | 0.4 | 4.0 | 4.0 | 0.001681 | 0.001674 | 99.6 |

A heuristic argument for the justification of the $z$-statistic is made in Appendix D.

Under the SW imputation, the relative bias is larger for smaller samples and for smaller response rates. The new imputation produces much smaller relative bias for the variance estimator. For large sample sizes both imputation methods produce negligible relative biases of the variance estimators.

TABLE 2
*Relative bias (RB) and the z-statistic of Rubin's variance estimators under the two different imputation schemes (50,000 samples)*

| Parameter | $n$ | $r/n$ | RB | | $z$-statistic | |
|---|---|---|---|---|---|---|
| | | | SW | New | SW | New |
| Mean | 20 | 0.8 | 0.0624 | 0.0008 | 9.33 | 0.12 |
| | | 0.6 | 0.1520 | −0.0045 | 21.27 | −0.66 |
| | | 0.4 | 0.3221 | 0.0074 | 37.23 | 0.97 |
| | 200 | 0.8 | −0.0086 | −0.0124 | −1.36 | −1.94 |
| | | 0.6 | 0.0040 | −0.0045 | 0.61 | −0.69 |
| | | 0.4 | 0.0155 | −0.0037 | 2.29 | −0.55 |
| Slope | 20 | 0.8 | 0.0706 | 0.0095 | 10.39 | 1.41 |
| | | 0.6 | 0.1560 | −0.0064 | 21.21 | −0.90 |
| | | 0.4 | 0.3418 | 0.0046 | 37.81 | 0.59 |
| | 200 | 0.8 | 0.0129 | 0.0107 | 2.03 | 1.67 |
| | | 0.6 | 0.0175 | 0.0031 | 2.69 | 0.48 |
| | | 0.4 | 0.0240 | 0.0047 | 3.60 | 0.70 |

TABLE 3
*Mean length and the coverage of 95% confidence intervals under the two different imputation schemes (50,000 samples)*

| Parameter | $n$ | $r/n$ | Mean length | | Coverage (%) | |
|---|---|---|---|---|---|---|
| | | | SW | New | SW | New |
| Mean | 20 | 0.8 | 1.1397 | 1.0996 | 95.4 | 95.0 |
| | | 0.6 | 1.5305 | 1.4032 | 95.9 | 94.7 |
| | | 0.4 | 2.2635 | 1.9213 | 96.4 | 94.7 |
| | 200 | 0.8 | 0.3232 | 0.3223 | 95.0 | 94.9 |
| | | 0.6 | 0.3949 | 0.3928 | 94.8 | 94.7 |
| | | 0.4 | 0.5235 | 0.5170 | 94.8 | 94.6 |
| Slope | 20 | 0.8 | 0.4169 | 0.4020 | 95.5 | 95.0 |
| | | 0.6 | 0.5646 | 0.5175 | 95.8 | 94.7 |
| | | 0.4 | 0.7897 | 0.6696 | 96.7 | 94.9 |
| | 200 | 0.8 | 0.1124 | 0.1121 | 95.2 | 95.1 |
| | | 0.6 | 0.1374 | 0.1364 | 95.0 | 95.0 |
| | | 0.4 | 0.1811 | 0.1789 | 95.0 | 94.7 |

Table 3 displays the mean lengths and the coverages of 95% confidence intervals. The confidence intervals are $(\hat{\theta} - t\sqrt{\hat{V}}, \hat{\theta} + t\sqrt{\hat{V}})$, where $t = t_{0.025,\nu}$ and $\nu$ is computed using the method of Barnard and Rubin (1999). The coverages of the confidence intervals are all close to the nominal level. For small sample sizes the confidence intervals based on the new imputation are slightly narrower than the confidence intervals based on the SW method. Interval estimation shows less dramatic results for small sample size than variance estimation. This is partly because the distributions of point estimates are bell-shaped and partly because the degrees of freedom of Barnard and Rubin (1999) attenuate the effect of small sample bias of the variance estimator.

**6. Discussion.** We study the mean and the covariance structure of the data set generated by the conventional multiple-imputation method under the regression model. Using the mean and the covariance structure of the multiply-imputed data set, we investigate the exact bias of Rubin's variance estimator. The bias of Rubin's variance estimator is negligible for large sample sizes, as discussed by SW, but the bias may be sizable for small sample sizes. We propose a simple modified imputation method that is more efficient than the SW method and makes Rubin's variance estimator unbiased. When applied to a general class of linear estimators, the proposed method produces more efficient estimates and has smaller bias for variance estimation than that of the SW method. In a simulation study we found that the bias of Rubin's variance estimator under the SW method is remarkably large for small sample sizes. The bias of Rubin's variance estimator under the new method is reasonably small in the simulation.

In practice, the small sample bias of the multiple imputation variance estimator is of special concern when the scale parameter $\sigma$ is generated with small degrees of freedom. One such example is stratified sampling, where the sample selection is performed independently across the strata. In a stratified sample the assumption of equal $\sigma$ across the strata is not a reasonable assumption. Thus, the scale parameters have to be generated independently within each stratum, using only the respondents in the stratum, which often makes the degrees of freedom very small even for a large data set. The new method will significantly reduce the bias in this case.

A commonly used imputation model is the cell mean model, where the study variables are assumed homogeneous within each cell. Under the cell mean model Rubin and Schenker (1986) considered various multiple imputation methods. Since the imputation is performed separately within each cell, the scale parameters are generated independently within each cell. Thus, the methods considered by Rubin and Schenker (1986) are subject to small sample biases. The biases can be significant when there are a small number of respondents within each cell. Recently, Kim (2002) proposed an alternative imputation method of making the variance estimator unbiased under the cell mean model.

## APPENDIX A

PROOF OF LEMMA 2.1.
By (3),

$$E(\mathbf{x}'_j \boldsymbol{\beta}^*_{(k)} | \mathbf{y}_r) = \mathbf{x}'_j \hat{\boldsymbol{\beta}}_r. \tag{A.1}$$

So

$$\mathrm{Cov}\,(Y_i, Y^{**}_{j(k)}) = \mathrm{Cov}\,\{Y_i, E(Y^{**}_{j(k)}|\mathbf{y}_r)\} = \mathrm{Cov}(Y_i, \mathbf{x}'_j \hat{\boldsymbol{\beta}}_r) = \mathbf{x}'_i (X'_r X_r)^{-1} \mathbf{x}_j \sigma^2.$$

Now, by (2),

$$E(\sigma^{*2}_{(k)}) = E(\sigma^{*2}_{(k)}|\mathbf{y}_r) = E(\lambda \hat{\sigma}^2_r) = \lambda \sigma^2, \tag{A.2}$$

where $\lambda = (r-p-2)^{-1}(r-p)$. So

$$\begin{aligned}
\mathrm{Cov}\,(\mathbf{x}'_i \boldsymbol{\beta}^*_{(k)}, \mathbf{x}'_j \boldsymbol{\beta}^*_{(k)}|\mathbf{y}_r) &= \mathrm{Cov}\,[E(\mathbf{x}'_i \boldsymbol{\beta}^*_{(k)}|\mathbf{y}_r, \sigma^*_{(k)}), E(\mathbf{x}'_j \boldsymbol{\beta}^*_{(k)}|\mathbf{y}_r, \sigma^*_{(k)})|\mathbf{y}_r] \\
&\quad + E[\,\mathrm{Cov}\,(\mathbf{x}'_i \boldsymbol{\beta}^*_{(k)}, \mathbf{x}'_j \boldsymbol{\beta}^*_{(k)}|\mathbf{y}_r, \sigma^*_{(k)})|\mathbf{y}_r] \\
&= E[\mathbf{x}'_i (X'_r X_r)^{-1} \mathbf{x}_j \sigma^{*2}_{(k)}|\mathbf{y}_r] \\
&= \mathbf{x}'_i (X'_r X_r)^{-1} \mathbf{x}_j \hat{\sigma}^2_r \lambda,
\end{aligned} \tag{A.3}$$

and, for $k \neq s$,

$$\mathrm{Cov}\,(\mathbf{x}'_i \boldsymbol{\beta}^*_{(k)}, \mathbf{x}'_j \boldsymbol{\beta}^*_{(s)}|\mathbf{y}_r) = 0. \tag{A.4}$$

Hence, by (A.1) and (A.3),

$$
\begin{aligned}
\text{Cov}(\mathbf{x}'_i\boldsymbol{\beta}^*_{(k)}, \mathbf{x}'_j\boldsymbol{\beta}^*_{(k)}) &= \text{Cov}\left[E(\mathbf{x}'_i\boldsymbol{\beta}^*_{(k)}|\mathbf{y}_r), E(\mathbf{x}'_j\boldsymbol{\beta}^*_{(k)}|\mathbf{y}_r)\right] \\
&\quad + E[\text{Cov}(\mathbf{x}'_i\boldsymbol{\beta}^*_{(k)}, \mathbf{x}'_j\boldsymbol{\beta}^*_{(k)}|\mathbf{y}_r)] \\
&= \text{Cov}(\mathbf{x}'_i\hat{\boldsymbol{\beta}}_r, \mathbf{x}'_j\hat{\boldsymbol{\beta}}_r) + E\{\mathbf{x}'_i(X'_rX_r)^{-1}\mathbf{x}_j\hat{\sigma}^2_r\lambda\} \\
&= (1+\lambda)\mathbf{x}'_i(X'_rX_r)^{-1}\mathbf{x}_j\sigma^2
\end{aligned}
$$
(A.5)

and, for $k \neq s$, by (A.1) and (A.4),

$$
\begin{aligned}
\text{Cov}(\mathbf{x}'_i\boldsymbol{\beta}^*_{(k)}, \mathbf{x}'_j\boldsymbol{\beta}^*_{(s)}) &= E[\text{Cov}(\mathbf{x}'_i\boldsymbol{\beta}^*_{(k)}, \mathbf{x}'_j\boldsymbol{\beta}^*_{(s)}|\mathbf{y}_r)] \\
&\quad + \text{Cov}\left[E(\mathbf{x}'_i\boldsymbol{\beta}^*_{(k)}|\mathbf{y}_r), E(\mathbf{x}'_j\boldsymbol{\beta}^*_{(s)}|\mathbf{y}_r)\right] \\
&= \text{Cov}(\mathbf{x}'_i\hat{\boldsymbol{\beta}}_r, \mathbf{x}'_j\hat{\boldsymbol{\beta}}_r) = \mathbf{x}'_i(X'X)^{-1}\mathbf{x}_j\sigma^2.
\end{aligned}
$$
(A.6)

By (4) we have, for $i \neq j$,

$$\text{Cov}(Y^{**}_{i(k)}, Y^{**}_{j(s)}) = \text{Cov}(\mathbf{x}'_i\boldsymbol{\beta}^*_{(k)}, \mathbf{x}'_j\boldsymbol{\beta}^*_{(s)}) \tag{A.7}$$

and

$$\text{Cov}(Y^{**}_{i(k)}, Y^{**}_{j(k)}) = \begin{cases} \text{Var}(\mathbf{x}'_i\boldsymbol{\beta}^*_{(k)}) + E(\sigma^{*2}_{(k)}), & \text{if } i = j, \\ \text{Cov}(\mathbf{x}'_i\boldsymbol{\beta}^*_{(k)}, \mathbf{x}'_j\boldsymbol{\beta}^*_{(k)}), & \text{if } i \neq j. \end{cases} \tag{A.8}$$

Therefore, inserting (A.2), (A.5) and (A.6) into (A.7) and (A.8), result (12) follows. □

## APPENDIX B

PROOF OF THEOREM 2.1. Before we calculate the variance of the multiple imputation estimator we provide the following matrix identity.

LEMMA B.1. *Let $X_n$ be an $n \times p$ matrix of the form $X'_n = (X'_r, X'_{n-r})$, where $X_r$ is an $r \times p$ matrix and $X_{n-r}$ is an $(n-r) \times p$ matrix. Assume that $X'_nX_n$ and $X'_rX_r$ are nonsingular. Then*

$$
\begin{aligned}
(X'_rX_r)^{-1} &= (X'_nX_n)^{-1} + (X'_nX_n)^{-1}X'_{n-r}X_{n-r}(X'_nX_n)^{-1} \\
&\quad + (X'_nX_n)^{-1}X'_{n-r}X_{n-r}(X'_rX_r)^{-1}X'_{n-r}X_{n-r}(X'_nX_n)^{-1}.
\end{aligned}
$$
(B.1)

PROOF. Using the identity [e.g., Searle (1982), page 261]

$$(D - CA^{-1}B)^{-1} = D^{-1} + D^{-1}C(A - BD^{-1}C)^{-1}BD^{-1} \tag{B.2}$$

with $A = I$, $B = X_{n-r}$, $C = X'_{n-r}$ and $D = X'_nX_n$, we have

$$
\begin{aligned}
(X'_rX_r)^{-1} &= (X'_nX_n)^{-1} + (X'_nX_n)^{-1}X'_{n-r}[I - X_{n-r}(X'_nX_n)^{-1}X'_{n-r}]^{-1} \\
&\quad \times X_{n-r}(X'_nX_n)^{-1}.
\end{aligned}
$$
(B.3)

Using (B.2) again with $A = X_n'X_n$, $B = X_{n-r}'$, $C = X_{n-r}$ and $D = I$, we have

$$\text{(B.4)} \quad [I - X_{n-r}(X_n'X_n)^{-1}X_{n-r}']^{-1} = I + X_{n-r}(X_r'X_r)^{-1}X_{n-r}'.$$

Inserting (B.4) into (B.3), we have (B.1). □

Note that $\hat{\boldsymbol{\beta}}_{M,n} = M^{-1}\sum_{k=1}^{M}\hat{\boldsymbol{\beta}}_{I(k),n}$ and the $\hat{\boldsymbol{\beta}}_{I(1),n}, \hat{\boldsymbol{\beta}}_{I(2),n}, \ldots, \hat{\boldsymbol{\beta}}_{I(M),n}$ are identically distributed. Thus,

$$\text{(B.5)} \quad \text{Var}(\hat{\boldsymbol{\beta}}_{M,n}) = (1 - M^{-1})\text{Cov}\,(\hat{\boldsymbol{\beta}}_{I(1),n}, \hat{\boldsymbol{\beta}}_{I(2),n}) + M^{-1}\text{Var}\,(\hat{\boldsymbol{\beta}}_{I(1),n}).$$

Define $\mathbf{a}_i = (X_n'X_n)^{-1}\mathbf{x}_i$ and $h_{ij} = \mathbf{x}_i'(X_r'X_r)^{-1}\mathbf{x}_j$. By (11) and (12)

$$\text{(B.6)} \quad \text{Cov}\,(\hat{\boldsymbol{\beta}}_{I(1),n}, \hat{\boldsymbol{\beta}}_{I(2),n}) = \sum_{i=1}^{r}\mathbf{a}_i\mathbf{a}_i'\sigma^2 + 2\sum_{i=1}^{r}\sum_{j=r+1}^{n}\mathbf{a}_i h_{ij}\mathbf{a}_j'\sigma^2 + \sum_{i=r+1}^{n}\sum_{j=r+1}^{n}\mathbf{a}_i h_{ij}\mathbf{a}_j'\sigma^2$$

and

$$\text{(B.7)} \quad \text{Var}\,(\hat{\boldsymbol{\beta}}_{I(1),n}) = \sum_{i=1}^{r}\mathbf{a}_i\mathbf{a}_i'\sigma^2 + 2\sum_{i=1}^{r}\sum_{j=r+1}^{n}\mathbf{a}_i h_{ij}\mathbf{a}_j'\sigma^2 + (1+\lambda)\sum_{i=r+1}^{n}\sum_{j=r+1}^{n}\mathbf{a}_i h_{ij}\mathbf{a}_j'\sigma^2 + \lambda\sum_{i=r+1}^{n}\mathbf{a}_i\mathbf{a}_i'\sigma^2.$$

By the definition of $\mathbf{a}_i$ and $h_{ij}$, we have

$$\text{(B.8)} \quad \sum_{i=1}^{r}\mathbf{a}_i\mathbf{a}_i' = (X_n'X_n)^{-1}X_r'X_r(X_n'X_n)^{-1},$$

$$\text{(B.9)} \quad \sum_{i=1}^{r}\sum_{j=r+1}^{n}\mathbf{a}_i h_{ij}\mathbf{a}_j' = (X_n'X_n)^{-1}X_{n-r}'X_{n-r}(X_n'X_n)^{-1} = \sum_{i=r+1}^{n}\sum_{j=r+1}^{n}\mathbf{a}_i\mathbf{a}_j'$$

and

$$\text{(B.10)} \quad \sum_{i=r+1}^{n}\sum_{j=r+1}^{n}\mathbf{a}_i h_{ij}\mathbf{a}_j' = (X_n'X_n)^{-1}X_{n-r}'X_{n-r}(X_r'X_r)^{-1}X_{n-r}'X_{n-r}(X_n'X_n)^{-1}.$$

Hence, inserting (B.8)–(B.10) into (B.6) and (B.7), and applying $X'_r X_r + X'_{n-r} X_{n-r} = X'_n X_n$ and (B.1), we have

$$\text{(B.11)} \qquad \text{Cov}\,(\hat{\boldsymbol{\beta}}_{I(1),n}, \hat{\boldsymbol{\beta}}_{I(2),n}) = (X'_r X_r)^{-1} \sigma^2$$

and

$$\text{(B.12)} \quad \text{Var}\,(\hat{\boldsymbol{\beta}}_{I(1),n}) = (X'_r X_r)^{-1} \sigma^2 + \lambda[(X'_r X_r)^{-1} - (X'_n X_n)^{-1}]\sigma^2.$$

Thus, (15) is proved by (B.5).

To show (17), because the $\hat{\boldsymbol{\beta}}_{I(1),n}, \hat{\boldsymbol{\beta}}_{I(2),n}, \ldots, \hat{\boldsymbol{\beta}}_{I(M),n}$ are identically distributed,

$$\text{(B.13)} \qquad E(B_{M,n}) = \text{Var}\,(\hat{\boldsymbol{\beta}}_{I(1),n}) - \text{Cov}\,(\hat{\boldsymbol{\beta}}_{I(1),n} \hat{\boldsymbol{\beta}}_{I(2),n}).$$

Thus, (17) is proved by inserting (B.11) and (B.12) into (B.13).

To show (16), we define $\tilde{\mathbf{Y}}'_{(k)} = (\mathbf{y}'_r, \mathbf{y}^{**'}_{(k)})$ to be the vector of the augmented data set at the $k$th repeated imputation, where $\mathbf{y}^{**}_{(k)} = (Y^{**}_{r+1(k)}, Y^{**}_{r+2(k)}, \ldots, Y^{**}_{n(k)})'$, $k = 1, 2, \ldots, M$. Then

$$(n-p)\hat{\sigma}^2_{I(k),n} = \tilde{\mathbf{Y}}'_{(k)}[I_n - X_n(X'_n X_n)^{-1} X'_n]\tilde{\mathbf{Y}}_{(k)}$$

and, under the regression model (1),

$$E\{\tilde{\mathbf{Y}}'_{(k)}[I_n - X_n(X'_n X_n)^{-1} X'_n]\tilde{\mathbf{Y}}_{(k)}\} = \text{trace}\,\{[I_n - X_n(X'_n X_n)^{-1} X'_n]V(\tilde{\mathbf{Y}}'_{(k)})\}.$$

By (11) and (12) we have

$$V(\tilde{\mathbf{Y}}'_{(k)}) = \begin{pmatrix} I_r & X_r(X'_r X_r)^{-1} X'_{n-r} \\ X_{n-r}(X'_r X_r)^{-1} X'_r & (1+\lambda)X_{n-r}(X'_r X_r)^{-1} X'_{n-r} + \lambda I_{n-r} \end{pmatrix} \sigma^2.$$

Let $P_x = X_n(X'_n X_n)^{-1} X'_n$ and $Q = V(\tilde{\mathbf{Y}}'_{(k)})\sigma^{-2} - I_n$. Then

$$\text{(B.14)} \quad \begin{aligned} E\{(n-p)\hat{\sigma}^2_{I(k),n}\} &= E\{(I - P_x)(Q_x + I_n)\sigma^2\} \\ &= \text{trace}\{I_n - P_x\}\sigma^2 + \text{trace}\{[I_n - P_x]Q\}\sigma^2. \end{aligned}$$

By the classical regression theory $\text{trace}\{I_n - P_x\} = n - p$. For the second term note that the left-upper $r \times r$ elements of $Q$ are all zeros. Define $C = (X'_n X_n)^{-1} X'_{n-r} X_{n-r}$ and $D = (X'_r X_r)^{-1} X'_{n-r} X_{n-r}$. Then

$$\text{(B.15)} \quad \begin{aligned} \text{trace}(Q) &= \text{trace}\,\{(1+\lambda)X_{n-r}(X'_r X_r)^{-1} X'_{n-r} + (\lambda-1)I_{n-r}\} \\ &= (1+\lambda)\,\text{trace}(D) + (\lambda-1)(n-r) \end{aligned}$$

and

$$\text{(B.16)} \quad \begin{aligned} \text{trace}(P_x Q) &= 2\,\text{trace}\,\{X_{n-r}(X'_n X_n)^{-1} X'_{n-r}\} \\ &\quad + (1+\lambda)\,\text{trace}\,\{X_{n-r}(X'_n X_n)^{-1} X'_{n-r} X_{n-r}(X'_r X_r)^{-1} X'_{n-r}\} \\ &\quad + (\lambda-1)\,\text{trace}\,\{X_{n-r}(X'_n X_n)^{-1} X'_{n-r}\} \\ &= (1+\lambda)\,\text{trace}(C + CD). \end{aligned}$$

Note that, using $X'_r X_r + X'_{n-r} X_{n-r} = X'_n X_n$, we have

$$\text{trace}(D - CD)$$
(B.17)
$$= \text{trace}\{[I - (X'_n X_n)^{-1} X'_{n-r} X_{n-r}](X'_r X_r)^{-1} X'_{n-r} X_{n-r}\}$$
$$= \text{trace}\{(X'_n X_n)^{-1} X'_{n-r} X_{n-r}\} = \text{trace}(C).$$

Thus, by applying (B.17) to (B.15) and (B.16), we have

(B.18) $$\text{trace}(Q - P_x Q) = (\lambda - 1)(n - r).$$

Therefore, inserting (B.18) into (B.14), we have (16). □

## APPENDIX C

PROOF OF THEOREM 4.1. The variance formula (29) directly follows by applying the multiple imputation variance formula (B.5) to $\hat{\theta}_{M,n}$ and using $\mathbf{a}_i = \alpha_i$ in (B.6) and (B.7).

To show (30), we decompose the total variance into three parts:

(C.1) $$\text{Var}(\hat{\theta}_{M,n}) = \text{Var}(\hat{\theta}_n) + \text{Var}(\hat{\theta}_{M,n} - \hat{\theta}_n) + 2\text{Cov}(\hat{\theta}_n, \hat{\theta}_{M,n} - \hat{\theta}_n).$$

To compute the bias of $W_{M,n}$ as an estimator of $\text{Var}(\hat{\theta}_n)$, we first express it as $\hat{V}_n = \mathbf{Y}'_n \Omega \mathbf{Y}_n$, where $\mathbf{Y}'_n = (\mathbf{y}'_r, \mathbf{y}'_{n-r})$ and $\Omega = [\Omega_{ij}]$; then the within-imputation variance term can be written $W_{M,n} = M^{-1} \sum_{k=1}^{M} \tilde{\mathbf{Y}}'_{(k)} \Omega \tilde{\mathbf{Y}}_{(k)}$. By $E(\tilde{\mathbf{Y}}_{(k)}) = E(\mathbf{Y}_n)$,

$$E(\hat{V}_{(k)}) = E(\tilde{\mathbf{Y}}'_{(k)}) \Omega E(\tilde{\mathbf{Y}}_{(k)}) + \text{trace}\{\Omega \text{Var}(\tilde{\mathbf{Y}}_{(k)})\}$$
$$= E(\hat{V}_n) + \text{trace}\{\Omega[\text{Var}(\tilde{\mathbf{Y}}_{(k)}) - \text{Var}(\mathbf{Y}_n)]\}.$$

By the unbiasedness of $\hat{V}_n$ and by the covariance structures in (11) and (12) we have

$$E(W_{M,n}) - \text{Var}(\hat{\theta}_n)$$
(C.2)
$$= 2 \sum_{i=1}^{r} \sum_{j=r+1}^{n} \Omega_{ij} h_{ij} \sigma^2 + (\sigma^2 + \lambda \sigma^2) \sum_{i=r+1}^{n} \sum_{j=r+1}^{n} \Omega_{ij} h_{ij}$$
$$+ (\lambda - 1) \sigma^2 \sum_{j=r+1}^{n} \Omega_{jj}.$$

For the $B_{M,n}$ term it can be shown that, using the same argument as for (B.6) and (B.7),

(C.3)
$$E(B_{M,n}) = \text{Var}(\hat{\theta}_{I(1),n}) - \text{Cov}(\hat{\theta}_{I(1),n} \hat{\theta}_{I(2),n})$$
$$= \left( \sum_{i=r+1}^{n} \sum_{j=r+1}^{n} \alpha_i \alpha_j h_{ij} + \sum_{j=r+1}^{n} \alpha_i^2 \right) \lambda \sigma^2.$$

By the covariance structures in (11)

$$\text{(C.4)} \qquad \text{Cov}(\hat{\theta}_{M,n}, \hat{\theta}_n) = \left( \sum_{i=1}^{r} \alpha_i^2 + \sum_{i=1}^{r} \sum_{j=r+1}^{n} \alpha_i \alpha_j h_{ij} \right) \sigma^2.$$

Thus, by (29), (C.4) and $\text{Var}(\hat{\theta}_n) = \sum_{i=1}^{n} \alpha_i^2 \sigma^2$,

$$\text{(C.5)} \qquad \begin{aligned} \text{Var}(\hat{\theta}_{M,n} - \hat{\theta}_n) &= \text{Var}(\hat{\theta}_{M,n}) + \text{Var}(\hat{\theta}_n) - 2\text{Cov}(\hat{\theta}_{M,n}, \hat{\theta}_n) \\ &= (1 + M^{-1}\lambda) \left( \sum_{i=r+1}^{n} \alpha_i^2 + \sum_{i=r+1}^{n} \sum_{j=r+1}^{n} \alpha_i \alpha_j h_{ij} \right) \sigma^2 \end{aligned}$$

and, by (C.3) and (C.5),

$$\text{(C.6)} \qquad \begin{aligned} &E[(1 + M^{-1}) B_{M,n}] - \text{Var}(\hat{\theta}_{M,n} - \hat{\theta}_n) \\ &= (\lambda - 1) \left[ \sum_{i=r+1}^{n} \alpha_i^2 + \sum_{i=r+1}^{n} \sum_{j=r+1}^{n} \alpha_i \alpha_j h_{ij} \right] \sigma^2. \end{aligned}$$

For the covariance term in (C.1) note that

$$\text{(C.7)} \qquad \begin{aligned} &\text{Cov}(\hat{\theta}_n, \hat{\theta}_{M,n} - \hat{\theta}_n) \\ &= \text{Cov}(\hat{\theta}_n, \hat{\theta}_{\infty,n} - \hat{\theta}_n) + \text{Cov}(\hat{\theta}_n, \hat{\theta}_{M,n} - \hat{\theta}_{\infty,n}) = 0, \end{aligned}$$

because the first term on the right-hand side of the above equality is zero by the congeniality condition (28) and the second term is also zero because

$$\text{Cov}(\hat{\theta}_n, \hat{\theta}_{M,n} - \hat{\theta}_{\infty,n}) = \text{Cov}(\hat{\theta}_n, \hat{\theta}_{M,n}) - \text{Cov}(\hat{\theta}_n, \hat{\theta}_{\infty,n}) = 0,$$

by the fact that $\hat{\theta}_{I(k)}$, $k = 1, 2, \ldots, M$, are identically distributed. Therefore, (30) follows from (C.2), (C.6) and (C.7). □

## APPENDIX D

**D.1. Justification for $z$-statistic in (35).** Let $(\hat{\theta}_i, \hat{V}_i)$, $i = 1, 2, \ldots, L$, be i.i.d. samples from a bivariate distribution $G(\hat{\theta}, \hat{V})$ with second moments. Then $E(\hat{V})$ is unbiasedly estimated by $E_L(\hat{V}) = L^{-1} \sum_{i=1}^{L} \hat{V}_i$ and $\text{Var}(\hat{\theta})$ is unbiasedly estimated by $(L-1)^{-1} L \times E_L[(\hat{\theta} - E_L(\hat{\theta}))^2] \doteq L^{-1} \sum_{i=1}^{L} (\hat{\theta} - E_L(\hat{\theta}))^2$, where $E_L(\hat{\theta}) = L^{-1} \sum_{i=1}^{L} \hat{\theta}_i$. Thus, by the central limit theorem,

$$\text{(D.1)} \qquad Z = \frac{E_L(\hat{V}) - E_L[(\hat{\theta} - E_L(\hat{\theta}))^2] - [E(\hat{V}) - \text{Var}(\hat{\theta})]}{\sqrt{\text{Var}\{E_L(\hat{V}) - E_L[(\hat{\theta} - E_L(\hat{\theta}))^2]\}}}$$

converges to a $N(0,1)$ distribution as $L \to \infty$. As $E_L(\hat{V}) - E_L[(\hat{\theta} - E_L(\hat{\theta}))^2] = L^{-1}\sum_{i=1}^{L}[\hat{V}_i - (\hat{\theta}_i - E_L(\hat{\theta}))^2]$, the variance term in the denominator of (D.1) is consistently estimated by

$$L^{-1}E_L\{[\hat{V} - (\hat{\theta} - E_L(\hat{\theta}))^2 - E_L[\hat{V} - (\hat{\theta} - E_L(\hat{\theta}))^2]]^2\}$$
$$\doteq L^{-1}E_L\{[\hat{V} - (\hat{\theta} - E_L(\hat{\theta}))^2 - E_L(\hat{V}) + \mathrm{Var}_L(\hat{\theta})]^2\}.$$

Thus, using Slutsky's theorem, the $z$-statistic in (35) converges to a $N(0,1)$ distribution under $H_0 : E(\hat{V}) = \mathrm{Var}(\hat{\theta})$ as $L \to \infty$.

**Acknowledgments.** The author thanks Wayne Fuller, Lou Rizzo and the referees for their useful comments and suggestions.

Department of Applied Statistics
Yonsei University
134 Sinchon-dong, Seodaemun-gu
Seoul 120-749
Korea
e-mail: kimj@yonsei.ac.kr